\documentclass[a4paper,10pt]{article}
\usepackage[utf8x]{inputenc}
\usepackage[russian]{babel}

\usepackage{amssymb}

\title{Специальные бор - зоммерфельдовы лагранжевы подмногообразия в алгебраических многообразиях.}
\author{Ник.А. Тюрин}

\begin{document}

\maketitle

\begin{abstract} Вводится новое понятие специальности лагранжевых подмногообразий, удовлетворяющих условию Бора - Зоммерфельда, для алгебраических
многообразий. Показывается, что такое условие позволяет строить конечномерные модули специальных бор - зоммерфельдвых лагранжевых подмногообразий относительно
любого обильного линейного расслоения. Конструкция может быть использована в исследованиях зеркальной симметрии.

\end{abstract}

Сущность зеркальной симметрии в наибольшей общности была выражена Ю.И. Маниным как ``двойственность симплектической геометрии и комплексной кэлеровой
геометрии'' (см. [1]). Два алгебраических кэлеровых многообразия $M, W$ понимаются как ``зеркальные партнеры'', если
некоторые производные объекты, получаемые в рамках алгебраической и симплектической геометрий $M, W$ перекрестно эквивалентны:
например, в Гомологической зекральной симметрии М. Концевича (см. [2]) производная категория когерентных пучков $M$ должна быть эквивалентна
категории Фукая - Флоера $W$ и наоборот.  

А.Н. Тюрин, посвятивший много лет исследованию стабильных векторных расслоений, преполагал более геометрическое соответствие:
соответствие между векторными расслоениями и лагранжевыми подмногообразиями (см., напр. [3]). А именно, для пары трифолдов $M, W$ четномерные когомологии
представляют классы Черна векторных расслоений, из которых составляются конечномерные многообразия модулей стабильных векторных расслоений,
а средние нечетные когомологии можно реализовать лагранжевыми подмногообразиями, из которых необходимо формировать конечномерные многообразия модулей,
и тогда появляется возможность сравнения этих двух серий многообразия модулей для определения той самой двойственности, которая составляет сущность зеркальной симметрии.
При этом главная проблема --- в ``бесконечности'' лагранжевой геометрии в противовес конечномерности геометрии алгебраической.

Эта проблема решается введением условий специальности на лагранжевы подмногообразия: развивая идею калиброванных лагранжевых циклов, Д. МакЛин и Н. Хитчин (см. [4])
предложили условие специальности на лагранжевы подмногообразия кэлеровых многообразий Калаби - Яу, которое приводило к конечномерным модулям. Сразу вслед за этим появилась SYZ - гипотеза для
трифолдов Калаби - Яу (см. [5]), объясняющая зеркальную симметрию в терминах слоений на специальные лагранжевы торы, однако многолетные попытки построить
специальные лагранжевы слоения на гладких многообразиях Калаби - Яу не увенчались успехом, что привело к постепенному спаду популярности SYZ - конструкции в кругах
математических физиков, сконцентрировавшихся на гомологическом подходе М. Концевича. Д. Ору в [6] пытался несколько оживить ``специальный'' подход к зеркальной симметрии,
обобщив понятие специального лагранжева подмногообразия на случай многообразий Фано, при этом специальность становится относительной --- условие специальности зависит
от выбора дивизора из антиканонической линейной системы на многообразии Фано (поэтому можно определить подход Ору как обобщение конструкции МакЛина и Хитчина на
многообразия Калаби - Яу ``с краем''). Для торической проективной плоскости исключение антиканонического дивизора, состоящего из трех прямых,
позволяет построить специальной слоение в смысле Ору, которое просто совпадает со стандартным слоением Лиувилля; Д. Ору построил специальное лагранжево слоение при
исключении приводимого дивизора --- прямой и коники --- и высказал гипотезу о существовании специального лагранжева слоения на дополнении к гладкой эллиптической кривой,
однако до сих пор эта гипотеза доказана не была. Кроме примера из [7], где было построено специальное в смысле Ору лагранжево слоение на многообразии флагов $F^3$,
больше ничего до сих пор в этом направлении сделано не было.

В работе [8] были развиты идеи лагранжева геометрического квантования, а именно --- были исследованы модули лагранжевых подмногообразий, удовлетворяющих
условию Бора - Зоммерфельда. По замечанию А.Н. Тюрина, условие Бора - Зоммерфельда ``трансверсально'' условию специальности в случае Калаби - Яу,
откуда можно было бы определять конечные инварианты для трифолдов Калаби - Яу, зеркальные инвариантам Кассонса, см. [3]. Мы развиваем эту идею,
вводя новое условие специальности на бор - зоммерфельдовы лагранжевы подмногообразия, что приводит к интересным наблюдениям.

Пусть $(M, \omega)$ --- компактное односвязное сиплектическое многообразие, удовлетворяющее условию Бора - Зоммерфельда: класс $[\omega] \in H^2(M, \mathbb{R})$
является целочисленным. Тогда зафиксируем данные предквантования: линейное расслоение $L \to M$ и эрмитову связность $a$ на нем, такие что
форма кривизны $F_a = 2 \pi i \omega$, при этом такая связность последним условием определена однозначно с точностью до каибровочного преобразования.
Лагранжево подмногообразие $S \subset M$ удовлетворяет условию Бора - Зоммерфельда (BS для краткости), если ограничение $(L, a)|_{S}$ допускает ковариантно постоянное сечение
$\sigma_S$. Пусть $s \in \Gamma(M, L)$ --- произвольное гладкое сечение $L$.

{\bf Определение.} {\it Назовем BS- лагранжево подмногообразие $S$ специальным относительно сечения $s$, если $s|_S$ нигде не обращается в нуль на $S$ и
коэффициент пропорциональности $\alpha(s, S)$, определяемый из равенства $s|_S = \alpha(S,s) \sigma_S$, имеет постоянный аргумент.}

Так как это определение не зависит от домножений сечения $s$ на любую ненулевую константу, то на самом деле мы определели соотношение, выделяющее ``цикл
инциденции'' в прямом произведении
$$
{\cal U}_{SBS} \subset \mathbb{P}(\Gamma(M, L)) \times {\cal B}_S,
$$
где последним символом обозначается многообразие модулей бор - зоммерфельдовых циклов фиксированного топологического типа, см. [8].
А именно, пара $(p, S)$ принадлежит ${\cal U}_{SBS}$,  если $S$ является специальным бор - зоммерфельдовым относительно сечения $s$, класс которого определен
точкой $p$ проективизации пространства сечений.  Естественным образом определены проекции $p_1, p_2$ на первый и второй прямые сомножители.

``Конечность'' множества специальных бор - зоммерфельдовых подмногообразий отражена в следующем факте:

{\bf Основная теорема.} {\it Слои проекции $p_1: {\cal U}_{SBS} \to \mathbb{P}(\Gamma(M, L))$ дискретны.}

Кроме того, нетрудно установить невырожденность дифференциала $d p_1$ на ${\cal U}_{SBS}$ и локальную сюрективность проеции $p_1$.
Отсюда следует

{\bf Утверждение.} {\it Пространство ${\cal U}_{SBS}$ обладает кэлеровой структурой, поднятой с $\mathbb{P}(\Gamma(M, L))$.}

Заметим, что в самом общем случае  это позволяет использовать  пространство ${\cal U}_{SBS}$  в обобщении квантования Сурьо - Костанта.

Однако нас интересует случай очень конкретный: пусть наше симплектическое многообразие $(M, \omega)$ обладает интегрируемой комплексной структурой
$I$, согласованной с $\omega$. Это значит, что $M$ есть алгебраическое многообразие с главной поляризацией, задаваемой голоморфным расслоением $L$.
Тогда мы получаем конечномерное подпространство $\mathbb{P}(H^0(M_I, L)) \subset \mathbb{P}(\Gamma(M, L))$ голоморфных сечений и редуцированный цикл инциденции
$$
{\cal M}_{SBS} \subset \mathbb{P}(H^0(M_I, L)) \times {\cal B}_S
$$
вместе с двумя естественными проекциями на прямые сомножители, которые мы снова обозначаем как $p_i$.

Тогда {\bf Основная теорема} может быть дополнена следующими утверждениями:

$\quad$ (1) образ $p_1({\cal M}_{SBS}) \subset \mathbb{P}(H^0(M_I, L))$ есть открытое подмножество в $\mathbb{P}(H^0(M_I, L))$;

$\quad$ (2) слои $p_1$ конечны;

$\quad$ (3) многообразие модулей ${\cal M}_{SBS}$ гладкое в общей точке многообразие, допускающее существование кэлеровой структуры.

{\bf Пример.} Возьмем в качестве $M$ простешее возможное компактное односвзяное симплектическое многообразие --- комплексную проективную прямую
$\mathbb{C} \mathbb{P}^1$ со стандартной комплексной структурой и кэлеровой формой, класс которой двойствен по Пуанкаре удвоенному классу точки
(это минимальный целочисленный класс, допускающий бор - зоммерфельдовы лагранжевы подмногообразия). Тогда расслоение $L$ изоморфно ${\cal O}(2)$,
любое голоморфное сечение с точностью до умножения на константу определяется парой точек (возможно совпадающих), в которых это сечение обращается в нуль.
Гладкая петля $\gamma \subset \mathbb{C} \mathbb{P}^1$ автоматически лагранжева из соображений размерности, и она удовлетворяет условию Бора - Зоммерфельда если и только если 
она делит $\mathbb{C} \mathbb{P}^1$ на две области равной симплектической площади. В этом случае ${\cal M}_{SBS}$ естественно изоморфно $\mathbb{C} \mathbb{P}^2 \backslash
Q$, где коника $Q$ есть образ $\mathbb{C} \mathbb{P}^1$ при вложении Веронезе. Более точно, в этом случае: (1) образом $p_1({\cal M}_{SBS}) \subset \mathbb{P}(H^0(M_I, L))
= \mathbb{C} \mathbb{P}^2$ будут сечения с некратными нулями; (2) слои $p_1$ состоят из единственных точек; (3) кэлерова структура с $\mathbb{C}\mathbb{P}^2$
поднимается до кэлеровой структуры на ${\cal M}_{SBS}$.  

Кроме того, на многообразии модулей ${\cal M}_{SBS}$ естественно определено  линейное расслоение ${\cal L}$ с эрмитовой связностью $A \in {\cal A}_h ({\cal L})$,
форма кривизны которой равна кэлеровой форме, поднятой с проективного пространства $\mathbb{P}(H^0(M_I, L))$, так что история только начинается...

Доказательство {\bf Основной теоремы}, детали теории специальных бор - зоммерфельдовых циклов и прочие примеры могут быть найдены в [9].

$$
$$

{\bf Литература:}

$$$$

[1] {\bf Ю.И. Манин}, {\it ``Предисловие ко третьему тому''}, Сборник избранных трудов А.Н. Тюрина в трех томах, Институт компьютерных технологий, Москва - Ижевск 2004; 

[2] {\bf M. Kontsevich}, {\it ``Homological algebra of mirror symmetry''}, Proceedings of ICM (Zurich, 1994), Birkhouser, Basel 1995, pp. 120 - 139;

[3] {\bf A.N. Tyurin}, {\it ``Geometric quantization and mirror symmetry''}, arXiv: math/9902027v1;

[4] {\bf N. Hitchin}, {\it ``Lectures on special lagrangian submanifolds''},  Winter school on Mirror Symmetry (Cambridge MA, 1999), AMS/IP Stud. Adv. Math,
23, AMS 2001, Providence,  pp. 151 - 182;

[5] {\bf A. Strominger, S.-T. Yau, E. Zaslow}, {\it ``Mirror symmetry is T - duality''}, Winter school on Mirror Symmetry (Cambridge MA, 1999), AMS/IP Stud. Adv. Math,
23, AMS 2001, Providence,  pp. 333 - 347;

[6] {\bf D. Auroux}, {\it ``Mirror symmetry and T - duality in the complement of an anticanonical divisor''}, J. Gokova Geom. Topol. (2007), pp. 51 - 91;

[7] {\bf Н.А. Тюрин}, {\it ``Специальные лагранжевы слоения многообразия флагов F3''}, Теор. Мат. Физ. 167: 2 (2011), стр. 193 - 205;

[8] {\bf А.Л. Городенцев, А.Н. Тюрин}, {\it ``Абелева лагранжева алгебраическая геометрия''}, Известия РАН, сер. мат. 65: 3 (2001), стр. 15 - 50;

[9] {\bf Н.А. Тюрин}, {\it ``Специальная геометрия лагранжевых циклов Бора - Зоммерфельда в алгебраических многообразиях''}, готовится к печати.

\end{document}